\documentclass{article}

\usepackage{arxiv}

\usepackage[utf8]{inputenc} 
\usepackage[T1]{fontenc}    
\usepackage{hyperref}       
\usepackage{url}            
\usepackage{booktabs}       
\usepackage{amsfonts}       
\usepackage{nicefrac}       
\usepackage{microtype}      
\usepackage{lipsum}		

\usepackage{tikz}                        
\usetikzlibrary{matrix,arrows}
\usepackage{amsmath}
\usepackage{amsthm}
\usepackage{amssymb}
\usepackage{mathtools}                  
\usepackage{dsfont}                     
\usepackage{amsxtra}                     
\usepackage{epsfig}
\usepackage{verbatim}
\usepackage{epigraph}                    

\theoremstyle{plain}
\newtheorem{thm}{Theorem} 
\newtheorem{prop}[thm]{Proposition}
\newtheorem{corollary}[thm]{Corollary}
\newtheorem{lemma}[thm]{Lemma}

\theoremstyle{definition}
\newtheorem{dfn}[thm]{Definition}

\newtheorem{example}[thm]{Example}

\theoremstyle{remark}
\newtheorem{remark}[thm]{Remark}
\newcommand{\bigslant}[2]{{\raisebox{.2em}{$#1$}\left/\raisebox{-.2em}{$#2$}\right.}}
\newcommand{\be}{\begin{equation}}
\newcommand{\ee}{\end{equation}}
\newcommand{\bdm}{\begin{displaymath}}
\newcommand{\edm}{\end{displaymath}}
\newcommand{\ba}[1]{\begin{array}{#1}}
\newcommand{\ea}{\end{array}}

\newcommand{\bcen}{\begin{center}}
\newcommand{\ecen}{\end{center}}
\newcommand{\btab}{\begin{tabular}}
\newcommand{\etab}{\end{tabular}}

\newcommand{\ra}{\rightarrow}
\newcommand{\lra}{\longrightarrow}

\newcommand{\Lra}{\Leftrightarrow}
\newcommand{\lmapsto}{\longmapsto}

\newcommand{\supp}{\ensuremath{\mathrm{supp}}}

\newcommand{\vertiii}[1]{{\left\vert\kern-0.25ex\left\vert\kern-0.25ex\left\vert #1 
		\right\vert\kern-0.25ex\right\vert\kern-0.25ex\right\vert}}
\newcommand{\PM}{\ensuremath{\mathcal{P}}} 
\newcommand{\B}{\ensuremath{\mathcal{B}}} 
\newcommand{\Sh}{\ensuremath{\mathcal{S}}} 

\newcommand{\R}{\ensuremath{\mathbb{R}}}
\newcommand{\N}{\ensuremath{\mathbb{N}}}

\title{Shape spaces in terms of Wasserstein geometry}

\date{October 22, 2025}	

\author{
	Bernadette Lessel \\
	University of Bonn\\
	\texttt{blessel@uni-bonn.de}
}

\begin{document}
\maketitle

\begin{abstract}
For a Polish space $X$, we define the Shape space $\Sh_p(X)$ to be the Wasserstein space $W_p(X)$ modulo the action of a subgroup $G$ of the isometry group $ISO(X)$ of $X$, where the action is given by the pushforward of measures. The Wasserstein distance can then naturally be transformed into a \emph{Shape distance} on Shape space if $X$ and the action of $G$ are proper. This is shown for example to be the case for complete connected Riemannian manifolds with $G$ being equipped with the compact-open topology. 
Before finally proposing a notion for tangent spaces on the Shape space $\Sh_2(\R^n)$, it is shown that $\Sh_p(X)$ is Polish as well in case $X$ and the action of $G$ are indeed proper. Also, the metric geodesics in $\Sh_p(X)$ are put in relation to the ones in $W_p(X)$.
\end{abstract}


\tableofcontents

\section{Introduction}

A Shape space in the sense of this article is a space whose elements are not defined by their position against an "absolute background space", as for example it is most trivially the case for $x\in\R^n$. Instead, an element of Shape space is defined by an intrinsic property that is subsequently called "shape". 
Whereas $x\in\R^n$ does not seem to have a (non-trivial) shape to begin with, the shape of a probability measure, on the other hand, can be regarded as being the the \emph{shape} of the spatial extension of its density. 
In this article, the formal approach to capture this idea is, however, indirect: the shape of a measure is thought to be what remains unchanged under the action of the isometry group of the underlying metric space. More concretely, a shape is going to be the equivalence class $[\mu]$ of a probability measure $\mu$ with respect to the action of that isometry group. The collection of all the equivalence classes $[\mu]$ will be called Shape space.

Elements in Shape space are put in relation to each other by the so-called Shape distance. This is, by definition, the smallest Wasserstein distance that can be realized between elements of two given equivalence classes of measures. Or, put differently, to measure the Shape distance of two measures $\mu$ and $\nu$, the smallest Wasserstein distance possible of $\mu$ to an element $g_\#\nu$ is taken, where $g$ is an element of the above mentioned isometry group and $g_\#\nu$ denotes the pushforward measure. 

This definition of Shape spaces, including the so-defined Shape distance, is motivated by the notion of Shape space from the physics literature. Shape spaces were therein pioneered by Julian Barbour and Bruno Bertotti in 1982 \cite{Barbour1982}, aiming at a dynamical description of the universe only through the relationships between objects, not relying on external reference frames or absolute time. The Shape dynamics program that subsequently followed builds on these relational principles by focusing on the intrinsic shapes formed by matter as the fundamental variables of the universe. A prominent role therein plays the homonymous Shape distance, a metric structure that gives the relative distance of shapes of configurations of matter. 

While picking up this idea of Shape spaces and the Shape distance from the context of Shape dynamics, the content of this article is, however, not any further related to this theory. Is is, in particular, mathematically completely different.

This article begins in the first section with recalling the necessary notions and properties of Wasserstein geometry. It does so following closely the introduction to the topic by \cite{User}. More details can be taken e.g. from \cite{Vilsmall} and \cite{Vilbig}. An extra subsection is dedicated to the tangent space-construction on Wasserstein spaces, in order to motivate the definition of tangent spaces on Shape space in section \ref{ch.infSh}.
Before the latter is being introduced, section \ref{ch.shape} introduces the notions of Shape spaces and Shape distances. An extra subsection is dedicated to Shape spaces on Riemannian manifolds. Subsequently, also geodesic and topological properties of Shape spaces are discussed. 

This preprint is a late spin-off from my dissertation "Shape space in terms of Wasserstein geometry and applications to quantum physics" at the University of Göttingen \cite{Lessel2018} in 2018, succeeding the preprint on "Differentiable maps between Wasserstein spaces" \cite{Lessel2020} from the same thesis.
\section{Wasserstein calculus}\label{ch.was}
%
Let $(X,d)$ be a Polish space and $\PM(X)$ the set of all probability measures on $X$ with respect to the Borel $\sigma$-algebra $\B(X)$. Instead of $(X,d)$ we will often just write $X$. A measurable map between two Polish spaces $T:X\rightarrow Y$ induces a map between the respective spaces of probability measures via the \emph{pushforward} $T_\#$ of measures: $T_\#:\PM(X)\rightarrow\PM(Y)$, $\mu\mapsto T_\#\mu$, where $T_\#\mu(A):=\mu(T^{-1}(A))$, for $A\in\B(Y)$. The \emph{support} of a measure $\mu$ is defined by $supp(\mu):=\{x\in X\mid \text{ every open neighbourhood of $x$ has positive $\mu$-measure}\}$. The Lebesgue measure on $\R^n$ is denoted by $\lambda$. 
%
\subsection{Wasserstein spaces $W_p(X)$}\label{subsec.Wass}
%
Wasserstein spaces are sets of probability measures together with the so called Wasserstein distance. To this end, we denote the set of probability measures which have finite $p$-th moment by $\PM_p(X)$, where $p\in[1,\infty)$:
$$\PM_p(X):=\{\mu\in\PM(X)\ \mid \int_X d^p(x_0,x)\ d\mu(x)<\infty,\ x_0\in X\}.$$

Furthermore, for $\mu, \nu$ in $\PM_p(X)$ we define
\begin{equation}\nonumber
Adm(\mu,\nu):=\lbrace \gamma\in\PM(X\times Y)\mid\pi^X_\#\gamma=\mu, \pi^Y_\#\gamma=\nu\rbrace,
\end{equation}
the so called \emph{admissible transport plans} between $\mu$ and $\nu$. Here, $\pi^X:X\times Y\rightarrow X$, $\pi^X(x,y)=x$. Similarly $\pi^Y$.

\begin{dfn}[\emph{\textbf{Wasserstein distances and Wasserstein spaces}}]
	Let $(X,d)$ be a Polish space and $p\in[0,\infty)$, then
	\begin{eqnarray*}\nonumber
		W_p: \PM_p(X)\times\PM_p(X)&\rightarrow& X\\\nonumber
		(\mu,\nu) &\mapsto& \left(\inf_{\gamma\in Adm(\mu,\nu)}\int_{X\times X}d^p(x,y)\ d\gamma(x,y)\right)^{1/p}
	\end{eqnarray*}
	is called the \emph{{$p$-th Wasserstein distance}}, or Wasserstein distance of order $p$. The tuple $(\PM_p(X),W_p)$ is called \emph{Wasserstein space} and is denoted by the symbol $W_p(X)$. 
\end{dfn}

The map $X\rightarrow W_p(X),\ x\mapsto\delta_x$ is an isometric embedding of the underlying Polish space $X$ into the Wasserstein space on $X$.

The fact that $W_p$ is indeed a metric distance, is a problem treated in optimal transport, where it is established that a minimizer for $$\inf_{\gamma\in Adm(\mu,\nu)}\int_{X\times X}d^p(x,y)\ d\gamma(x,y)$$ actually exists. Such a minimizer is called \emph{optimal transport plan}. In case a plan $\gamma\in Adm(\mu,\nu)$ is induced by a measurable map $T:X\rightarrow Y$, i.e. in case $\gamma=(Id,T)_\#\mu$, $T$ is called \emph{transport map}. Then, $T_\#\mu=\nu$.

One can show that $W_p(X)$ is compact, in case $X$ is. On the other hand, whenever $X$ is unbounded, then $W_p(X)$ cannot be locally compact. Furthermore, $W_p(X)$ is complete and separable and the Wasserstein distance $W_p$ metrizes the weak convergence in $\PM_p(X)$. In particular, $W_p$ is continuous on $\PM_p(X)\times \PM_p(X)$.

\begin{dfn}[\emph{\textbf{Weak convergence in $\PM_p(X)$}}]
	A sequence $(\mu_k)_{k\in\N}\subset\PM(X)$ is said to \emph{converge weakly} to $\mu\in\PM(X)$ if and only if $\int\varphi d\mu_k\ra\int\varphi d\mu$ for any bounded continuous function $\varphi$ on $X$. This is denoted by $\mu_k\rightharpoondown\mu$.
	A sequence $(\mu_k)_{k\in\N}\subset\PM_p(X)$ is said to converge weakly to $\mu\in\PM_p(X)$ if and only if for $x_0\in X$ it is:
	\begin{itemize}
		\item[1)] $\mu_k\rightharpoondown\mu$ and
		\item[2)] $\int d(x_0,x)^p d\mu_k(x)\ra\int d(x_0,x)^p d\mu(x)$.
	\end{itemize}
	This is denoted by $\mu_k\rightharpoonup\mu$.
\end{dfn}

Let us now briefly have a look at the geodesic structure of Wasserstein spaces. The notion of a geodesic in important for the construction of the tangent spaces and in Subsection \ref{subsec.geod} we will compare geodesics in Wasserstein space with geodesics in Shape space.

\begin{dfn}[\emph{\textbf{Constant speed geodesic}}]\label{df.geod}
	A curve $(\gamma_t)_{t\in [0,1]}$, $\gamma_0\neq\gamma_1$, in a metric space $(X,d)$ is called a \emph{constant speed geodesic} or \emph{metric geodesic} in case that
	\begin{equation}\label{eq.geod}
	d(\gamma_t,\gamma_s)=|t-s|d(\gamma_0,\gamma_1)\ \ \forall t,s\in [0,1]. 
	\end{equation}
\end{dfn}

Equivalently, Condition \ref{eq.geod} can be stated with an "$\leq$" instead of an "$=$". 

We will often abbreviate curves $(\gamma_t)_{t\in [0,1]}$ by omitting the brackets and the interval of definition, i.e we will often just write $\gamma_t$ instead.

Since $x\mapsto \delta_x$ is an isometry, for every constant speed geodesic $\gamma_t$ in $X$, $\delta_{\gamma_t}$ is a constant speed geodesic in $W_p(X)$.
\begin{dfn}[\emph{\textbf{Geodesic space}}]
	A metric space $(X,d)$ is called \emph{geodesic} if for every $x,y\in X$ there exists a constant speed geodesic $(\gamma_t)_{t\in [0,1]}$ with $\gamma_0=x$ and $\gamma_1=y$. 
\end{dfn}

$W_2(X)$ is geodesic whenever $(X,d)$ is. If $X$  is a Hilbert space, $\mu_t$ is a constant speed geodesic in $W_2(X)$ if and only if there exists an optimal transport plan $\gamma\in Adm(\mu_0,\mu_1)$ such that $\mu_t=((1-t)\pi^1+t\pi^2)_\#\gamma.$
If further $\gamma$ is induced by a map $T$, this latter formula simplifies to $\mu_t=((1-t)Id+tT)_\#\mu_0.$

Curves of the form $\mu_t=((1-t)\pi^1+t\pi^2)_\#\gamma$ are called \emph{displacement interpolation} from $\pi^1_\#\gamma$ to $\pi^2_\#\gamma$ through $\gamma$. The previous statement thus tells that within the Wasserstein structure, it is more natural to interpolate measures on the level of measurable sets than on the level of measures itself: The curve $\mu_t=(1-t)\mu_0+t\mu_1$, where interpolation is done by just shifting mass from one measure to the other, has infinite length in Wasserstein space and is, although being continuous, not absolutely continuous. In $((1-t)\pi^1+t\pi^2)_\#\gamma$ the prescription is that a set $A$ has the amount of mass that $\gamma$ is assigning to the set of all those points $(x,y)$ with interpolation $(1-t)x+ty$ being an element of $A$.

In the following Definition, let $Geod(X)$ be the set of all constant speed geodesic on $X$.

\begin{dfn}[\emph{\textbf{Non branching space}}]
	A metric space $(X,d)$ is called \emph{non branching}, if the following map is injective for every $t\in(0,1)$
	\begin{eqnarray}\nonumber
	Geod(X) &\longrightarrow& X\times X\\\nonumber
	\ \gamma\ &\longmapsto& (\gamma_0,\gamma_t).
	\end{eqnarray}
\end{dfn}

	If $(X,d)$ is a complete, separable, locally compact and non branching geodesic space, then also $W_p(X)$ is non branching for $p\in(0,\infty)$. The converse is also true: If $W_p(X)$ is non branching, then so is $X$. 
A proof of this can be found in \cite{Vilbig} (Corollary 7.32). According to \cite{User}, the local compactness condition on $X$ is, however, not necessary.
%
\subsection{Tangent spaces on $W_2(M)$}\label{sec.dyn}
%
Let now the Polish space $X$ be a smooth, connected and complete Riemannian manifold manifold $(M,h)$. We will often omit the subscript $2$ in $W_2(M)$ and write $W(M)$ instead.
Furthermore, we equip the set of measurable sections in $TM$, which we will denote by $\Gamma(TM)$, with an $L^2$-topology. That means, for $v\in \Gamma(TM)$ we define
\begin{equation*}
\|v\|_{L^2(\mu)}:=\sqrt{\int_M h(v,v)\ d\mu}
\end{equation*}
and  \begin{equation*}
L^2(M,\mu):=\bigslant{\{v\in \Gamma(TM)\mid \|v\|_{L^2(\mu)} <\infty\}}{\sim}.
\end{equation*}
Here, two vector fields are considered to be equivalent in case they differ only on a set of $\mu$-measure zero. $L^2(M,\mu)$ is a Hilbert space with the canonical scalar product. We will often write $L^2(\mu)$ if it is clear to which manifold $M$ it is referred to.

The (infinite dimensional) manifold structure that is commonly used on $W(M)$ is not a smooth structure in the sense of e.g. \cite{Kriegl_1997} where infinite dimensional manifolds are modeled on convenient vector spaces. The differentiable structure on $W(M)$, that will be introduced below, rather consists of ad hoc definitions accurately tailored to optimal transport and the Wasserstein metric structure which only mimic conventional differentiable and Riemannian behavior.

Instead of starting with a smooth manifold structure, on Wasserstein spaces one starts with the notion of a tangent space. 
Traditionally, the basic idea of a tangent vector at a given point is that it indicates the direction a (smooth) curve will be going infinitesimally from that point. Then, the set of all such vectors which can be found to be tangential to some curve at a given fixed point are collected in the tangent space at that point. On $W(M)$, however, there is no notion of smooth curves. But there is a notion of metric geodesics. In case the transport plan for the optimal transport between two measures is induced by a map $T$, the interpolating geodesic on Hilbert spaces can be written as $\mu_t=((1-t)Id+tT)_\#\mu_0$, thus being of the form $\mu_t=F_{t\#}\mu_0$. More generally, on Riemannian manifolds optimal transport between $\mu_0$ and $\mu_t$ can be achieved by $\mu_t=F_{t\#}\mu_0$, $F_t=\exp(t\nabla\varphi)$( see e.g \cite{Vilbig}, Chapter 12). In these cases, $F_t$ is injective and locally Lipschitz for $0<t<1$ (\cite{Vilsmall}, Subsubsection 5.4.1).
It is known from the theory of characteristics for partial differential equations that curves of this kind solve the weak continuity equation, together with the vector field to which integral lines $F_t$ corresponds. 
\begin{dfn}[\emph{\textbf{Continuity equation}}]\label{dfn.conti}
	Given a family of vector fields $(v_t)_{t\in [0,T]}$, a curve $\mu_t:[0,T]\rightarrow W(M)$ is said to solve the \emph{continuity equation}
	\begin{equation}\label{eq.cont}
	\partial_t\mu_t+\nabla\cdot(v_t\mu_t)=0
	\end{equation}
	in the weak sense, in case
	\begin{equation}\label{eq.wcont}
	\int_0^T\int_M\left(\frac{\partial}{\partial t}\varphi(x,t)+h(\nabla\varphi(x,t),v_t(x))\right)\ d\mu_t(x)dt =0
	\end{equation}
	holds true for all $\varphi\in C_c^\infty\left((0,T)\times M\right)$.
\end{dfn}
\begin{thm}\label{thm.conti}
	Let $(F_t)_{t\in[0,T)}$ be a family of maps on $M$ such that $F_{t}:M\rightarrow M$ is a bijection for every $t\in[0,T)$, $F_0=Id$ and both $(t,x)\mapsto F_t(x)$ and $(t,x)\mapsto F_t^{-1}(x)$ are locally Lipschitz on $[0,T)\times M$. Let further $v_t(x)$ be a family of velocity fields on $M$ such that its integral lines correspond to the trajectories $F_t$, and $\mu$ be a probability measure. Then $\mu_t=F_{t\#}\mu$ is the unique weak solution in $\mathcal{C}\left([0,T),\PM(M)\right)$ of $\frac{d}{dt}\mu_t+\nabla\cdot(v_t\mu_t)=0$ with initial condition $\mu_0=\mu$. Here, $\PM(M)$ is equipped with the weak topology.
\end{thm}
%
%
Theorem \ref{thm.conti} is taken from \cite{Vilsmall} where it is labeled as Theorem 5.34. 

It it possible to characterize the class of curves on $W(M)$ that admit a velocity in the manner of \ref{dfn.conti} (\cite{User}). 
\begin{dfn}[\emph{\textbf{Absolutely continuous curve}}]
	Let $(E,d)$ be an arbitrary metric space and $I$ an interval in $\R$. A function $\gamma:I\rightarrow E$ is called \emph{absolutely continuous, a.c.,} if there exists a function $f\in L^1(I)$ such that 
	\begin{equation}\label{eq.ac}
	d(\gamma(t),\gamma(s))\leq\int_t^s f(r)dr,\ \ \ \forall s,t\in I, t\leq s.
	\end{equation}
\end{dfn}
\begin{dfn}[\emph{\textbf{Metric derivative}}] \label{df.md}
	The \emph{metric derivative} $|\dot{\gamma}|(t)$ of a curve $\gamma:[0,1]\rightarrow E$ at $t\in (0,1) $ is given as the limit
	\begin{equation}
	|\dot{\gamma}|(t)=\lim_{h\rightarrow 0}\frac{d(\gamma(t+h),\gamma(t))}{|h|}.
	\end{equation}
\end{dfn}
Every metric geodesic is absolutely continuous and $|\dot{\gamma}|(t)=d(\gamma(0),\gamma(1))$.

 It is known that for absolutely continuous curves $\gamma$, the metric derivative exists for a.e. $t$. It is an element of $L^1(0,1)$ and, up to sets of zero Lebesgue-measure, the minimal function satisfying equation \eqref{eq.ac} for $\gamma$. In this sense absolutely continuous functions enable a generalization of the fundamental theorem of calculus to arbitrary metric spaces.

\begin{thm}[\emph{\textbf{Differential characterization of a.c. curves}}]\label{them.ac}
	Let $\mu_t:[0,1]\rightarrow W(M)$ be an a.c. curve. Then there exists a Borel family of vector fields $(v_t)_{t\in[0,1]}$ on $M$ such that the continuity equation \eqref{eq.wcont} holds and 
	\begin{equation*}
	\|v_t\|_{L^2(\mu_t)}\leq |\dot{\mu_t}| \text{  for a.e. } t\in(0,1).
	\end{equation*}
	Conversely, if a curve $\mu_t:[0,1]\rightarrow W(M)$ is such that there exists a Borel family of vector fields $(v_t)_{t\in[0,1]}$ with $\|v_t\|_{L^2(\mu_t)}\in L^1(0,1)$, together with which it satisfies \eqref{eq.wcont}, then there exists an a.c. curve $\tilde{\mu}_t$ being equal to $\mu_t$ for a.e. $t$ and satisfying
	\begin{equation*}
	|\dot{\tilde{\mu}}_t|\leq\|v_t\|_{L^2(\tilde{\mu}_t)} \text{  for a.e. } t\in(0,1).
	\end{equation*} 
	\qed
\end{thm}
 In the following we call a pair $(\mu_t,v_t)$ which together solves the continuity equation an \emph{a.c. couple}. We further want to call $v_t$ an \emph{accompanying vector field} for $\mu_t$.

%

As seen in Theorem \ref{them.ac}, every absolutely continuous curve in $W(M)$ admits an $L^1(dt)$ family of $L^2(\mu_t)$-vector fields together with which the continuity equation is satisfied. This family of vector fields $v_t$ can be seen as being tangential to the curve. However, $v_t$ is not unique, there are many vector fields which allow for the same motion of the density: Adding another family $w_t$ with the ($t$-independent) property $\nabla(w_t \mu_t)=0$ to $v_t$  does not alter the equation. Theorem \ref{them.ac} provides a natural criterion to choose a unique element among the $v_t's$. According to this theorem, there is at least one family $v_t$ such that $|\dot{\mu}_t|=\|v_t\|_{L^2(\mu_t)}$ for almost all $t$, i.e. that is of minimal norm for almost all $t$. Linearity of \eqref{eq.ac} with respect to $v_t$ and the strict convexity of the $L^2$-norms ensure the uniqueness of this choice, up to sets of zero measure with respect to $t$. We want to call such a couple $(\mu_t,v_t)$, where $v_t$ is the unique minimal tangent family for an a.c. curve $\mu_t$, a \emph{tangent couple}. 

It then seems reasonable to define the tangent space at point $\mu$ as the set of $v\in L^2(M,\mu)$ with $\|v\|_\mu\leq\|v+w\|_\mu$ for all $w\in L^2(M,\mu)$ such that $\nabla(w\mu)=0$. This condition for $v\in L^2(M,\mu)$, however, is equivalent to saying that $\int_M h(v,w)\ d\mu=0$ for all $w\in L^2(M,\mu)$ with $\nabla(w\mu)=0$. This in turn is equivalent to the following, which we will take as the definition of the tangent space.

\begin{dfn}[\emph{\textbf{Tangent space $T_\mu W(M)$}}] \label{dfn.tspace}
	The \emph{tangent space} $T_\mu W(M)$ at point $\mu\in W(M)$ is defined as
	\begin{equation}\label{eq.tangentspace}
	T_\mu W(M):=\overline{\{\nabla\varphi\mid \varphi\in \mathcal{C}^\infty_c(M)\}}^{L^2(M,\mu)}\subset L^2(M,\mu).
	\end{equation}	
\end{dfn}

\begin{remark}\label{cor.tan}
	It can be shown that if $(\mu_t,v_t)$ is an a.c. couple such that $\|v_t\|_{L^2(\mu_t)}\in L^1(0,1)$, then $(\mu_t, v_t)$ is a tangent couple if and only if $v_t\in T_{\mu_t}W(M)$ for almost every $t\in(0,1)$.
\end{remark}
It is not difficult to see that $\text{dim }T_{\delta} W(M)= \text{dim }M$, whereas in most of the cases $\text{dim }T_\mu W(M)=\infty$. However, every probability measure can be approximated by a sequence of measures with finite support (see \cite{Vilbig} Thm 6.18), so that in each neighborhood of every measure there is an element $\mu$ with $\text{dim }T_\mu W(M)<\infty$.
Furthermore, as long as $\mu$ is concentrated on an at most countable set, it is $T_{\mu} W(M)=L^2(M,\mu)$ (see \cite{Gigli}, Remark 1.33).

%

\section{Shape spaces $\Sh_p(X)$}\label{ch.shape}
Let $(X,d)$ be again a Polish metric space and $ISO(X)$ its isometry group. Further, $p\in\N_{\ge 1}$.

\subsection{Shape distances $D_p(\mu,\nu)$}\label{sec.shape}

Every subgroup $G$ of $ISO(X)$ acts naturally on $\PM_p(X)$ by means of the pushforward of measures:
	\begin{eqnarray*}\nonumber
		G\times \PM_p(X) &\rightarrow& \PM_p(X)\\\nonumber
		(g,\mu) &\mapsto& g\mu:=g_\#\mu.
	\end{eqnarray*}



\begin{dfn}[\emph{\textbf{Shape distance}}]
The \emph{p-th Shape distance} $D_p(\mu,\nu)$ of $\mu,\nu\in\PM_p(X)$ is defined as 
\begin{equation*}
D_p(\mu,\nu):= \inf_{g\in G}{W_p(g\mu,\nu)}. 
\end{equation*}
\end{dfn}

In the following we will show that this \emph{Shape distance} is a pseudometric on $\PM_p(X)$ modulo $G$. If $X$ is such that the closed balls are compact and the action of $G$ on $X$ is proper, $D_p(\cdot,\cdot)$ is really a metric on that space, as we will see afterwards.

We recall the definition of a pseudometric. 

\begin{dfn}[\emph{\textbf{Pseudometric}}]
	A function $d:X\times X\rightarrow\R_{\geq0}$ between a set $X$ is called a \emph{pseudometric} in case for all $x,y,z\in X$, $d(x,x)=0$, $d(x,y)=d(y,x)$ and $d(x,z)\leq d(x,y)+d(y,z)$. Other than for a metric, it is not demanded that $d(x,y)=0$ if and only if $x=y$.
\end{dfn}

\begin{dfn} \label{dfn.equiv}
Two elements in $\PM_p(X)$ are defined to be \emph{equivalent} if there exists an element in $G$ which maps them into each other:
\begin{equation*}
\mu\sim\nu :\Lra\ \exists\ g\in G: \mu=g\nu.
\end{equation*}
The quotient space with respect to this equivalence relation is denoted by $\Sh_p(X):=\bigslant{\PM_p(X)}{\sim}$.
\end{dfn}
\begin{dfn}[\emph{\textbf{Shape space}}]
The space $\Sh_p(X)$, together with $D_p(\cdot,\cdot)$ is called \emph{Shape space}. If $p=2$, we will often just write $\Sh(X)$ instead of $\Sh_2(X)$
\end{dfn}
\begin{remark}
It is evident that Definition \ref{dfn.equiv} indeed yields an equivalence relation on $\PM_p(X)$. Also, we observe 
that the following map is well defined, where we denote the equivalence class of $\mu$ by $[\mu]$:
\begin{eqnarray}\nonumber
D_p(\cdot,\cdot):\Sh_p(X)\times \Sh_p(X) &\lra& \R\\ \nonumber
                 \left([\mu],[\nu]\right) \;\;\quad &\lmapsto& D_p(\mu,\nu). 
\end{eqnarray}
For this reason, we always mean a map on $\Sh_p(X)\times \Sh_p(X)$ when we write expressions such 
as $D_p(\mu,\nu)$. We state it explicitly when we mean a map on $\PM_p(X)$.
\end{remark}

\begin{thm}\label{thm.pseudo}
The p-th Shape distance $D_p(\cdot,\cdot)$ is a pseudometric on $\Sh_p(X)$.
\end{thm}
To prove this theorem we formulate the following two lemma:

\begin{lemma}\label{bij}
For topological spaces $X,Y$ let $f,g:X\rightarrow Y$ be injective measurable maps such that $f(X)$ and $g(X)$ are closed. Also, let  $\widetilde{Adm}(f_\#\mu,g_\#\nu):=\{\Pi\in\PM\left(f(X)\times g(X)\right)\mid \pi^1_\#\Pi=f_\#\mu, \pi^2_\#\Pi=g_\#\nu)\}$. Then $(f\times g)_\#(\cdot)$ is a bijective map from $Adm(\mu,\nu)$ to $\widetilde{Adm}(f_{\#}\mu,g_{\#}\nu)$ with inverse map $(f^{-1}\times g^{-1})_\#(\cdot)$. 
\end{lemma}

\begin{proof}
First we need to see that the maps are well-defined. For this let $\Pi\in Adm(\mu,\nu)$. It is $\supp\left((f \times g)_\# \Pi\right)\subset f(X)\times g(X)$, since the images of $X$ under $g$ and $f$ are closed, which means that we can naturally identify $(f \times g)_\# \Pi$ as an element of $\PM(f(X)\times g(X))$. Furthermore it is
\begin{equation*}
(\pi^1_\#(f \times g)_\# \Pi)(A)=\Pi(f^{-1}(A)\times X)=\mu(f^{-1}(A))=f_{\#}\mu(A),
\end{equation*}
for $A\in\B(f(X))$. A Similar calculation holds true for the projection onto the second component and for $(f_1^{-1}\times f_2^{-1})_\#(\cdot)$. The latter map is well-defined since every map is surjective on its image.
For being inverse we observe that for $A,B\in\B(X)$
\begin{equation*}
\left( (f_1^{-1}\times f_2^{-1})_\#(f_1\times f_2)_\#\Pi\right) (A,B)=\left( (f_1\times f_2)_\#\Pi\right) (f_1(A)\times f_2(B))=\Pi(A,B)
\end{equation*}
and the other way round.
\end{proof}

\begin{remark}
For every $\Pi\in Adm(f_\#\mu, g_\#\nu)$ it is $\supp\ \Pi\subset f(X)\times g(X)$ and in particular
\begin{eqnarray}\nonumber
W_p(f_\#\mu,g_\#\nu)^p &=& \inf_{\Pi\in Adm(f_\#\mu,g_\#\nu)}\int_{Y\times Y}d^p(x,y)\ d\Pi(x,y)\\ \nonumber
                  &=&\inf_{\Pi\in \widetilde{Adm}(f_\#\mu,g_\#\nu)}\int_{f(X)\times g(X)}d^p(x,y)\ d\Pi(x,y).
\end{eqnarray}
Of course, $\widetilde{Adm}(f_\#\mu,g_\#\nu)=Adm(f_\#\mu,g_\#\nu)$ in case $f$ and $g$ are also surjective.
\end{remark}

\begin{lemma}\label{sym}
For all $g\in G$ it is
\begin{equation*}
W_p(g\mu,\nu)=W_p(\mu,g^{-1}\nu).
\end{equation*}
\end{lemma}

\begin{proof}
According to Lemma $\ref{bij}$ the map $(g\times g)_\#(\cdot)$ gives a bijective correspondence of the sets $Adm(g\mu,\nu)$ and $Adm(\mu,g^{-1}\nu)$. We can use this to write 
\begin{eqnarray}\nonumber
W_p(g\mu,\nu)^p&=&\inf_{\Pi\in Adm(g\mu,\nu)}{\int_{X}d(x,y)^p d\Pi(x,y)}\\ \nonumber
             &=&\inf_{\bar{\Pi}\in Adm(\mu,g^{-1}\nu)}{\int_X d(x,y)^p d[(g\times g)_\#\bar{\Pi}](x,y)}\\ \nonumber
             &=&\inf_{\bar{\Pi}\in Adm(\mu,g^{-1}\nu)}{\int_X d(x,y)^p d\bar{\Pi}(x,y)}=W_p^p(\mu,g^{-1}\nu).
\end{eqnarray}
\end{proof}
\begin{corollary}\label{cor.isub}
From Lemma $\ref{sym}$ it follows that $G$ acts isometric also on the $p$-th Wassersteinspace: 
$W_p(g\mu,g\nu)=W_p(\mu,g^{-1}g\nu)=W_p(\mu,\nu)$. $\Box$
\end{corollary}
\begin{corollary}
Let $G_\#:=\{g_\#:W(M)\rightarrow W(M)\mid g\in G\}$. It is a subgroup of $ISO(W(M))$ and group-isomorphic to $G$.
\end{corollary}
\begin{proof}
$G_\#\subset ISO(W(M))$ follows from \ref{cor.isub}. Subgroup properties are immediate. The group homomorphism between $G$ and $G_\#$ is given by $\varphi:G\rightarrow G_\#;\ g\mapsto g_\#$. It is surjective by definition and injective due to the following reasoning:
Assume that $g_1,g_2\in G$ are such that $g_{1\#}=g_{2\#}$. Then for every $x \in M$ we have 
\[ \delta_{g_1(x)} = g_{1\#}\delta_x = g_{2\#}\delta_x = \delta_{g_1(x)},\] which means that $g_1(x) = g_2(x)$.
\end{proof}
\begin{proof}[Proof of Theorem \ref{thm.pseudo}]
Let $[\mu],[\nu]$ and $[\sigma]$ be elements of $\Sh_p(X)$. 
Of course $D_p(\mu,\nu)=0$ if $[\mu]=[\nu]$ with optimal $g=1_G$ and $D_p(\mu,\nu)\geq0$ in any other case, by definition.
Then we show that $D$ is symmetric:
Since $W_p$ is symmetric and from Lemma $\ref{sym}$ we know that 
\begin{equation*}
D_p(\mu,\nu)=\inf_{g\in G}W_p(g\mu,\nu)=\inf_{g\in G}W_p(\mu,g\nu)=\inf_{g\in G}W_p(g\nu,\mu)=D_p(\nu,\mu).
\end{equation*}
For the validity of the triangle inequality we see that for all $\sigma\in W_p(M)$ and all $g_1,g_2\in G$,
\begin{eqnarray}\nonumber
D_p(\mu,\nu)&=&\inf_{g\in G}W_p(g\mu,\nu)\leq W_p(g_2g_1\mu,\nu))= W_p(g_1\mu,g_2^{-1}\nu) \\ \nonumber
          &\leq& W_p(g_1\mu,\sigma)+W_p(\sigma,g_2^{-1}\nu)= W_p(g_1\mu,\sigma)+W_p(g_2\sigma,\nu).
\end{eqnarray}
Therefore, 
\begin{eqnarray}\nonumber
D_p(\mu,\nu)=\inf_{g\in G}W_p(g\mu,\nu)&\leq&\inf_{g_1,g_2\in G}\left(W_p(g_1\mu,\sigma)+W_p(g_2\sigma,\nu)\right) \\ \nonumber
                                       &=& \inf_{g_1\in G}W_p(g_1\mu,\sigma)+\inf_{g_2\in G}W_p(g_2\sigma,\nu)\\ \nonumber
                                       &=& D_p(\mu,\sigma)+D_p(\sigma,\nu).
\end{eqnarray}
\end{proof}
Actually, we have proven something far more general. The only properties we have used for Theorem $\ref{thm.pseudo}$ is that the group action of $g^{-1}$ behaves as a sort of adjoint action to the action of $g$ with respect to the distance on the original metric space and that the formula for the prospective distance on the quotient space is given as an infimum over all possible distances elements of two orbits can have. We will formulate this finding in the following corollary.
\begin{corollary}\label{cor.gen}
Let the metric space $(X,d)$ and $G$ be such that $d(gx,y)=d(x,g^{-1}y)$ for every $g\in G$ and $D(x,y)=\inf_{g\in G}d(gx,y)$. Then $D$ is a pseudometric on $X/G$.
\end{corollary}

\begin{remark}\label{rem.rema}
The condition $d(gx,y)=d(x,g^{-1}y)$ is necessary for $D$ to be well-defined on $X/G$ and lets us freely put the $g$ at one of the two arguments.
\end{remark}
\begin{remark}
Note that it is not clear under which conditions $D(\cdot,\cdot)$ would be symmetric if $G$ would not be a group.
\end{remark}
\begin{remark}
$D(\cdot,\cdot)$ is not a metric in this generality, because a minimizer for $\inf_{g\in G}d(\cdot,\cdot)$ need not exist. But as we will see, it is possible to show this minimizing property for proper actions. 
\end{remark}
\begin{remark}
One might wonder if one can bring $W_p(\mu,g\nu)$ and $W_p(\mu,\nu)$ in any algebraic relationship to each other, for every $\mu,\nu\in\PM(X)$ and fixed $g$. However, a simple homogeneous relationship can not be established, as one can see at the example $X$ being the real numbers, $p=1$ and $gx:=x+2.$ If $\mu=\delta_x$ and $\nu=\delta_y$, then $W(\mu,\nu)=d(x,y).$ Choose here for example $x=3$ and $y=0$, then $3\cdot W(\mu,g\nu)=W(\mu,\nu).$ If, on the other hand, $x=2$ then $0\cdot W(\mu,\nu)=W(\mu,g\nu).$
\end{remark}
\subsection{Existence of a minimizer for the Shape distance-problem}\label{s.minimizer}
The goal of this section is to prove the existence of a minimizer for the Shape distance-problem. The basic idea is that the problem should be reduced to having a continuous map on a compact subset of the isometry group and the minimizer of this map should be exactly the minimizer of the shape distance-problem. This is why we start with the following lemma and its important corollary.
\begin{lemma}\label{cont}
Let $G$ be a topological group with continuous group action on $X$. Then the group action of $G$ on $W_p(X)$ is continuous as well. 
\end{lemma}
\begin{proof}
Let $\left( (g_i),(\mu_i)\right)_{i\in I}$ be a net in $G\times W_p(X)$ with $\left( g_i,\mu_i\right)\rightarrow(g,\mu)$ (see Appendix \ref{A.nets}). As we know from Section $\ref{sec.was}$, the Wasserstein distance metrizes weak convergence, i.e. we have to show that $W_p(g_i\mu_i,g\mu)\ra0$ as $\left( g_i,\mu_i\right)\rightarrow(g,\mu)$. For this we see that
\begin{equation*}
W_p(g_i\mu_i,g\mu)=W_p(\mu_i,g_i^{-1}g\mu)\leq W_p(\mu_i,\mu)+W_p(\mu,g_i^{-1}g\mu)\quad\forall i\in I.
\end{equation*}
We know already that $\mu_i\rightharpoonup\mu$. But it is also true that $g_i^{-1}g\mu\rightharpoonup\mu$: 
Since $g_i\rightarrow g$ in $G$ and $G$ acts continuously on $X$, the functions $g^{-i}g(x)$ converge pointwise against $Id_G(x)$. Using dominated convergence theorem yields $g_i^{-1}g\mu\rightharpoonup\mu$ and thus the continuity of the Wasserstein distance on $W_p(X)$.
\end{proof}

\begin{remark}
In case $G\subset ISO(X)$, there is always at least one topology on $G$ such that the action on X is continuous, namely the compact-open topology. Generally, the following holds: Is $X$ a locally compact and $Y$ any topological space and is further $H\subset C(X,Y)$, where $C(X,Y)$ is the set of continuous functions $f:X\rightarrow Y$. Then one can show (\cite{Schubert}, p.74) that the compact-open topology is the coarsest topology on $H$ such that the map $H\times X\rightarrow Y,\ (\Phi,x)\mapsto \Phi(x)$ is continuous. 
\end{remark}


\begin{corollary}\label{cor.cont}
Let $G$ and $(X,d)$ be as in Lemma $\ref{cont}$. Then the map $G\ra\R$;\:$g\mapsto W_p(g\mu,\nu)$ is continuous for every choice of $\mu,\nu\in W_p(X)$.
\end{corollary}
\begin{proof}
We know from the literature that $W_p(\cdot,\cdot)$ is continuous on $W_p(X)\times W_p(X)$. With this, $W_p(\cdot,\nu)$ is 
continuous on $W_p(X)$ for every choice of $\nu$. Since we know from Lemma $\ref{cont}$ that the group action of $G$ on $W_p(X)$ is also continuous, we have continuity as a composition of continuous functions.
\end{proof}

The map considered in Corollary $\ref{cor.cont}$ is of course the map for which we want to find the existence of a minimizer. The complement of the compact set on which this map is supposed to be considered for this should consist of all those elements in $G$ which transports the two given measures ``too far away'', so that the distance between them is increasing instead of decreasing. 
For this, it is useful to consider spaces on which most of the mass of probability measures is concentrated on compact balls, to have better control over their Wasserstein distances.  To synthesize this, the action should ``take care'' of compact sets appropriately. The property we will use for this is \emph{proper}.
\begin{dfn}[\emph{\textbf{Proper action}}]\label{dfn.proper}
Let $X$ be a metric space. A continuous action $G\times X\rightarrow X$ of a group $G$ on $X$ is called \emph{proper}, if the map $G\times X\rightarrow X\times X,\ (g,x)\mapsto (gx,x)$ is proper.
\end{dfn}
%
%
\begin{remark}
	The definition for an action to be proper is equivalent to saying that for every $K_1, K_2\subset X$ compact the set $\{g\in G\mid gK_1\cap K_2\neq\emptyset\}$ is compact in $G$.	
	
	One can show that if $G$ acts proper on $X$ then $X/G$ is a Hausdorff space.
\end{remark}
\begin{dfn}[\emph{\textbf{Proper metric space}}]\label{dfn.properspace}
A metric space is called \emph{proper} if every bounded closed set is compact.
\end{dfn}
\begin{remark}\label{rem.rem}
	\begin{itemize}
		\item 	A metric space $(X,d)$ being proper is equivalent to $x\mapsto d(x_0,x)$ being a proper map for any $x_0\in X$.
		\item   A proper metric space is automatically locally compact (which is why $W_p(X)$ cannot be proper in general). This fact enables us to use, the for our situation suitable, Definition \ref{dfn.proper} of proper actions on locally compact spaces.
		\item  Any proper space is complete.
		\item  One can view ``proper'' as a finiteness condition, recalling that finite dimensional Hilbert spaces are always proper, but in infinite dimensions one can show that the closed unit ball is not sequentially compact.
	\end{itemize}
\end{remark}
Now we have all the ingredients we need to formulate the theorem.

\begin{thm}[\emph{\textbf{Existence of a minimizer}}]\label{min}
Let the Polish metric space $(X,d)$ be proper. Additionally let the topology and the action of $G$ be such that it acts properly on $X$. In this case a minimizer for the problem $\displaystyle inf_{g\in G}W_p(g\mu,\nu)$ exists for every choice of $\mu,\nu\in\PM_p(X)$.
\end{thm}
\begin{corollary}\label{metric}
In the setting of Theorem $\ref{min}$ the Shape distance $D_p(\cdot,\cdot)$ is a metric.
\end{corollary}
\begin{proof} [Proof of Corollary \ref{metric}]
Since we already know that $D_p(\cdot,\cdot)$ is a pseudometric, we only have to make sure that if $D_p(\mu,\nu)=0$ then automatically 
$[\mu]=[\nu]$. From Theorem \ref{min} we know that there exists an element $g_{min}\in G$ such that $0=D_p(\mu,\nu)=W_p(g_{min}\mu,\nu)$. Since $W_p(\cdot,\cdot)$ is a metric we can conclude that $g_{min}\mu=\nu$ and thus $[\mu]=[\nu]$.
\end{proof}
\begin{remark}
If there exists no minimizer of $inf_{g\in G}W_p(g\mu,\nu)$ in general, one can always convert the pseudometric $D_p(\cdot,\cdot)$ into a metric via the so called metric identification. For this procedure one identifies points which have zero pseudodistance with each other. 
\end{remark}
We will use the following definition in the proof for the existence of a minimizer.
\begin{dfn}
Let $X$ be a topological space. A sequence of subsets $(K_n)_{n\in\N}\subset X$ is called \emph{exhaustion of X} if $\bigcup_{n\in\N}K_n=X$ and $K_{n}\subset K_{n+1}$.
\end{dfn}
\begin{proof}[Proof of Theorem \ref{min}]
If $X$ is bounded, then it is automatically compact and then also $G$ is compact.
So from now on, let $X$ be unbounded.

According to the prerequisites, for every $R\in\R,\ x\in X:\ \overline{B_R(x)}$ is compact and we have an exhaustion of $X$ by balls centered around an arbitrary element $x\in X$.
Now let $\mu,\nu\in \PM_p(X)$ and $x\in X$. Using the Monotone Convergence Theorem, for every $\epsilon>0$ there exists an element $R\in\R$ such that $\mu(B_R(x))> 1-\epsilon$ and $\nu(B_R(x))> 1-\epsilon$. Then let $C:=W_p(\mu,\nu)^p$ and choose $g\in G$ such that $g \overline{B_{R'}(x)}\cap \overline{B_{R'}(x)}=\emptyset$, where $\displaystyle R'=kR,\ k>\frac{1}{2R}\sqrt[p]{\frac{C}{1-2\epsilon}}+1$. To fulfill this, we have to choose $\epsilon<1/2$.


If no such $g$ exists, we know automatically that $G$ is compact and the proof ends here. 
Otherwise, for $B':=B_{R'}(x)$, it is
\begin{eqnarray}\nonumber
& & W_p(g\mu,\nu)^p\\\nonumber 
&=& \inf_{\Pi\in Adm(g\mu,\nu)}\left(\int_{gB\times B} d^p(x,y)d\Pi(x,y)+\int_{X\times X\setminus gB\times B} d^p(x,y) d\Pi(x,y)\right)\\ \nonumber
                &\geq& (2(R'-R))^p \inf_{\Pi\in Adm(g\mu,\nu)}\Pi(gB\times B)\geq (2(k-1)R)^p(1-2\epsilon)\\ \nonumber
                &>&  C
\end{eqnarray}

This means that the minimizer has to be found in the subset $\{g\in G\mid g \overline{B_{R'}(x)}\cap \overline{B_{R'}(x)}\neq \emptyset\}\neq\emptyset$. Since all closed balls of $(X,d)$ are compact and we assumed the group action to be proper, this subset is compact, so the existence of a minimizer is guaranteed.

\end{proof}
\begin{remark}
One might have the idea that, instead of demanding that there exists an exhaustion of $X$ by compact balls, it is enough to say that there exists a compact exhaustion, i.e. to say that $X$ is $\sigma$-compact. But in this case, it is not clear how to control the distances of elements of the set $gK_{n'}$ to elements of the set $K_{n'}$, where $K_n$ are the elements of the exhaustion.
\end{remark}
To confirm intuition, we additionally give the following proposition, which is implicitly used in the preceding proof:
\begin{prop}
Let $\mu\in\PM_p(X)$ and $g\in G$. Then $\supp(g\mu)=g\ \supp(\mu)$. In particular, two measures cannot be equivalent to each other if each support is not an isometric image of the respective other support.
\end{prop}
\begin{proof}
Let $g\in G,\ \mu\in\PM_p(X)$. 
To establish this proposition, we need the following two identities: 
\begin{itemize}
\item[1.)] $\{A\in \mathcal{B}(X)\mid \mu(g^{-1}(A^c))=0\}=\{gA\in\mathcal{B}(X)\mid \mu(A^c)=0\}$.
\item[2.)] Let $I$ be an index family, we then have $\displaystyle\bigcap_{i\in I}gA_i=g\bigcap_{i\in I}A_i$ for $A_i\subset X$.
\end{itemize}
Thereby is $g A:=\{ga\mid a\in A\}$, for $A\subset X$. Proof of the first identity:\\
Let $A\in \mathcal{B}(X)$, then $A=gg^{-1}A$, since $g$ is bijective. Then we have $\mu((g^{-1}A)^c)=\mu(g^{-1}A^c)=0$. For the other way round, let $gA\in\mathcal{B}(X)$, then $\mu(g^{-1}(gA)^c)=\mu(A^c)=0$. We will use the second identity without further justification. Now we can see that 
\begin{eqnarray}\nonumber
\supp\ g\mu &=& \overline{\bigcap_{g\mu(A^c)=0}A}=\overline{\bigcap_{\mu(g^{-1}A^c)=0}A}=\overline{\bigcap_{\mu(A^c)=0}gA} \\\nonumber
            &=& \overline{g\bigcap_{\mu(A^c)=0}A} =g\ \overline{\bigcap_{\mu(A^c)=0}A}=g\ \supp\mu
\end{eqnarray}
\end{proof}
\subsection{Shape spaces $\Sh_p(M)$ on Riemannian manifolds}\label{sec.shRiem}
For applications in physical scenarios, it is important to know whether Theorem \ref{min} is applicable to the Euclidean space, or more generally to Riemannian manifolds. 

In the realm of Riemannian geometry, one a priori distinguishes between two notions of isometry. Since every connected Riemannian manifold $(M,h)$ is a metric space, being equipped with the Riemannian distance, there is on the one hand the notion of isometry in the sense defined in the beginning of this chapter. 
On the other hand, a diffeomorphism $f:(M_1,h_1)\rightarrow (M_2,h_2)$ is called an isometry, in case it satisfies $h_1=f^*h_2$, where $^*$ denotes the pullback of tensors. Fortunately, every such diffeomorphism in this latter sense is also an isometry in the metric sense (whenever $(M_1,h_1)=(M_2,h_2)$). The theorem by Myers and Steenrod guarantees that the other way round is true as well \cite{Myers_1939}.

We further remind that Riemannian manifolds are separable and locally compact metric spaces by definition and the metric topology coincides with the manifold topology. However, in general Riemannian manifolds are not complete, i.e. Polish, but for our purposes it is of course necessary to require this from now on.

There are two further things to check. First, whether every complete Riemannian manifold is proper and secondly, whether the action of the isometry group is proper. Regarding the first issue we can make use of the following Theorem by Hopf and Rinow:

\begin{thm}[\emph{\textbf{Hopf-Rinow theorem, \cite{Rinow_1964}}}]\label{thm.HopfRinow}
Let $(M, h)$ be a Riemannian manifold. Then the following is equivalent:
\begin{itemize}
\item[i)] $M$ is a complete metric space.
\item[ii)] Every closed bounded set is compact.
\end{itemize}
\end{thm}

Now we have to have a look at the isometry group of $(M,h)$ and its action on $M$.
\begin{thm}[\cite{DantzigWaerden}] 
Let $(M,h)$ be a connected Riemannian manifold, then its isometry group $G$, equipped with the compact-open topology, is a locally compact topological group with continuous group action on $M$. In case $M$ is compact, so is $G$.
\end{thm}

A proof of this can also be found in \cite{KobNom}.

We now state the next significant theorem.
\begin{thm}
Let $(M,h)$ be a connected Riemannian manifold and let its isometry group $G$ be equipped with the compact-open topology. Then the action of $G$ on $M$ is proper.
\end{thm}
A proof of this can be found in \cite{Schlecht}, p. 17, Theorem 2.22.
\begin{corollary}
In case the Polish space $X$ is a complete connected Riemannian manifold and its isometry group is equipped with the compact-open topology, a minimizer for the Shape distance-problem exists and $D_p(\cdot,\cdot)$ is a metric distance on $\Sh_p(X)$.
\end{corollary}

The Euclidean space is of course a complete and connected Riemannian manifold.

\begin{remark}\begin{itemize}
\item According to \cite{Schlecht}, in case $M$ is the union of more than one connected components, the action of the isometry group need not to be proper in general.
\item The isometry group of a connected Riemannian manifold $M$ is a closed subset of the group of diffeomorphisms of $M$ with respect to the compact-open topology (see \cite{KobNom} p. 48, or \cite{Schlecht}).
\end{itemize}
\begin{remark}
Infinite dimensional (separable) Banach spaces don't serve as an example for our situation since closed unit balls are not compact in general. In this sense, one can regard the condition of the compactness of closed balls in Theorem \ref{min} as a sort of finiteness condition (compare also Remark \ref{rem.rem}). 
\end{remark}
\end{remark}
\subsection{Some topological properties of $\Sh_p(X)$}

\begin{thm}\label{thm.topology} 
Let the metric space $(X,d)$ and the subgroup of the isometry group $G$ be such that $d(gx,y)=d(x,g^{-1}y)$ for every $g\in G$. Let furthermore $D(x,y):=\inf_{g\in G}d(gx,y)$ be a metric on $X/G$ in such a way that the infimum is always realized by an element $g\in G$. Then for $\displaystyle\pi:X\rightarrow X/G,x \mapsto [x]$ it is equivalent:
\begin{itemize}
\item[1.)] $A$ is open in $(X/G,D)$
\item[2.)]$\pi^{-1}(A)$ is open in $(X,d)$. 
\end{itemize}
\end{thm}
\begin{remark}
Compare the conditions in this theorem (and the following) to the ones in \ref{cor.gen} and in \ref{rem.cond}.
\end{remark}
\begin{corollary}
The topology on Shape space coming from the Shape distance is the final topology with respect to the projection map $\pi$ (i.e. the quotient topology). In particular, $\pi$ is continuous. \qed
\end{corollary}
\begin{proof}[Proof of Theorem \ref{thm.topology}]
Since in metric spaces the open $r$-Balls are a basis of the topology, we only have to prove \ref{thm.topology} for $A=B^{X/G}_r([x]):=\{[y]\in X/G\mid D([x],[y])<r)\}.$ Then we have for $B^X_r(x):=\{y\in X\mid d(x,y)<r\}$,
\begin{eqnarray}\nonumber
\pi^{-1}\left(B^{X/G}_r([x])\right) &=& \{y\in X\mid [y]\in B^{X/G}_r([x])\}\\\nonumber
                                      &=& \{y\in X\mid\exists\ g\in G:d(gx,y)<r\}\\\nonumber
                                      &=& \bigcup_{g\in G}\{y\in X\mid d(gx,y)<r\}\\\nonumber
                                      &=& \bigcup_{g\in G} B_r^X(gx).
\end{eqnarray}
For the other direction let $A\subset X/G$ such that $\pi^{-1}(A)=\bigcup_{i\in I}B_{r_i}^X(x_i).$ Then 
$$A=\pi\left(\bigcup_{i\in I}B_{r_i}^X(x_i)\right)=\bigcup_{i\in I}\pi(B_{r_i}^X(x_i)).$$
So the question reduces to whether $\pi(B_{r}^X(x))\in\mathcal{O}_{X/G}$ for an arbitrary $x\in X$ and $r\in\R$:
\begin{eqnarray}\nonumber
\pi(B_{r}^X(x)) &=& \pi\left(\{y\in X\mid d(x,y)<r\}\right)\\\nonumber
                    &=& \{\pi(y)\in X/G\mid d(x,y)<r\}\\\nonumber
                    &=& \{[y]\in X/G\mid \exists\ \tilde{y}\in [y]: d(x,\tilde{y})<r\}\\\nonumber
                    &=& \{[y]\in X/G\mid D([x],[y])<r)\}\\\nonumber
                    &=& B_{r}^{X/G}(x)
\end{eqnarray}
\end{proof}
\begin{thm}\label{thm.complete}
Let the situation be as in Theorem \ref{thm.topology}. Then $(X/G,D)$ is complete if $X$ is complete.
\end{thm}
\begin{proof}
Let $([x_n])_{n\in\N}\subset X/G$ be a Cauchy sequence. The idea of the proof is to define a Cauchy sequence in $X$ such that its induced sequence in the space of equivalence classes is a subsequence of the original Cauchy sequence $([x_n])_{n\in\N}$. This is of use, since a Cauchy sequence converges if it has a converging subsequence. 

Now, for $\epsilon_0=(1/2)^0$ there is an $N_0\in\N$ such that for every $m,n>N_0$, $D([x_n],[x_m])< \epsilon_0$. Define $y_0$ to be an arbitrary element of $[x_{m_0}]$, for an $m_0\ge N_0$. Further, for $\epsilon_1=(1/2)^1$ there is an $N_1\in\N$ such that for every $m,n>N_1\ge m_0$, $D([x_n],[x_m])< \epsilon_1$. Define $y_1$ to be an element of $[x_{m_1}]$, for an $m_1\ge N_1$ such that $d(y_0,y_1)=D([y_0],[y_1])$, which is smaller $1$ in particular, according to the choice of $y_0$ and $y_1$. 
Successively, the sequence $(y_n)\subset X$ is chosen according to this scheme. In particular, for every $m\in\N$ it is 
$ d(y_0,y_m)<\sum_{i=0}^{\infty}{\epsilon_i}<\sum_{i=0}^{\infty}{(1/2)^i}<\infty$. Since $\epsilon_i>0 \forall i$ and $\lim_{n\rightarrow\infty}\epsilon_i=0$, we have 
$\lim_{n \rightarrow \infty}\sum_{i=n}^{\infty}{\epsilon_i}=0$ monotoneously and therefore for
every $\epsilon>0$ there is an $N\in\N$ such that for every $m,n>N$, $d(y_n,y_m)<\epsilon$. Since $X$ is complete, there is an $y\in X$ such that $\lim_{n\rightarrow\infty}y_n=y$.
Finally, $\lim_{n\rightarrow\infty}D([y_n],[y])\leq \lim_{n\rightarrow\infty}d(y_n,y)=0$.
\end{proof}
\begin{remark}\label{rem.cond}
The condition $d(gx,y)=d(x,g^{-1}y)$ is necessary for $D$ to be well-defined on $X/G$ and lets us freely put $g$ in one of the two arguments. And the existence of a minimum for $D(x,y)=\inf_{g\in G}d(gx,y)$ is needed in the proof to ensure the existence of the chosen elements $y_i,\ i>0$.
\end{remark}
For the next, resulting, corollary we will additionally use that the $p-$th Wasserstein space $W_p(X)$ is complete, whenever $(X,d)$ is (which we assume always to be the case). 
\begin{corollary}
$(\Sh_p(X),D_p)$ is complete. \qed 
\end{corollary}
\begin{prop}
Let the situation be like above. Then, if $(X,d)$ is separable then $(X/G,D)$ is separable as well.
\end{prop}

\begin{proof}
Let $A\subset X$ be countable and dense. Then $[A]:=\{[a]\in X/G \mid a\in A\}$ is also countable. Regarding denseness, let $M\in X/G$. Then for every $x\in M$ there exists a sequence $(x_n)_{n\in\N}\subset A$ such that $d(x_n,x)\longrightarrow 0$. But then also $D([x_n],[x])\longrightarrow 0$.
\end{proof}

\begin{corollary} \label{cor.polSh}
$\Sh_p(X)$ is a Polish space.\qed
\end{corollary}

In the following, by $(X,d_X)\cong(Y,d_Y)$ we mean that the metric spaces $(X,d_X)$ and $(Y,d_Y)$ are isometric, i.e. there is a bijective isometry between these spaces. Also, we denote the isometry group of $(X,d_X)$ by $G^X$ and the Shape Distance on $W_p(X)/G^X$ by $D_p^{G^X}$. $G^Y$ and $D_p^{G^Y}$ analogously.

\begin{lemma}\label{lem.iso}
	Let $(X,d_X)$ an $(Y,d_Y)$ be Polish spaces. If $I:(X,d_X)\hookrightarrow(Y,d_Y)$ is a distance-preserving map, i.e. $d_Y(I(x),I(y)=d_X(x,y)$ for all $x,y\in X$, then 
	$I_\#:W_p(X)\hookrightarrow W_p(Y)$ is a distance-preserving map as well.
\end{lemma}

\begin{proof}
We are going to use \ref{bij}, with that same notation. For this, it is necessary that $X$ is complete such that $I(X)$ is closed in $Y$
\begin{eqnarray} \nonumber
W_p^Y(I_\#\mu,I_\#\nu)^p &=& \inf_{\pi\in Adm(I_\#\mu,I_\#\nu)}\int_{Y\times Y}{d_Y^p(x,y)d\pi(x,y)}\\\nonumber
                       &=& \inf_{\pi\in \widetilde{Adm}(I_\#\mu,I_\#\nu)}\int_{I(X)\times I(X)}{d_Y^p(x,y)d\pi(x,y)}\\ \nonumber
                       &=&  \inf_{\pi\in Adm(\mu,\nu)}\int_{I(X)\times I(X)}{d_Y^p(x,y)dI_\#\pi(x,y)}\\ \nonumber
				       &=& \inf_{\pi\in Adm(\mu,\nu)}\int_{X\times X}{d_X^p(x,y)d\pi(x,y)}\\ \nonumber
	                   &=& W_p^X(\mu,\nu).
\end{eqnarray}
\end{proof}

\begin{corollary}
From $(X,d_X)\cong(Y,d_Y)$ follows $W_p(X)\cong W_p(Y)$.
\end{corollary}

\begin{proof}
Every distance preserving map $I$ is injective, surjectivity of $I_\#$ follows from $I$ being bijective: Let $\mu\in\PM(Y)$, then $I^{-1}_\#\mu$ is mapped to $\mu$ by $I_\#$.
\end{proof}

\begin{prop} \label{prop.iso}
If $I:(X,d_X)\hookrightarrow(Y,d_Y)$ is a distance-preserving map such that $I(X)$ is closed, then $[I_\#]:(W_p(X)/G^X, D^{G^X}_p)\hookrightarrow (W_p(Y)/G^Y, D^{G^Y}_p),\ [\mu]\mapsto [I_\#\mu]$ is a distance-preserving map as well. 
\end{prop}

\begin{proof}
We start by showing that $[I_\#]$ is well-defined. For this we first note that the map $G^X\rightarrow G^{I(X)}, g\mapsto I\circ g\circ I^{-1}$ is a well-defined bijective map. $G^{I(X)}$ is meant to be the full isometry group of the metric space $I(X)$, which can be identified with a subgroup of $G^Y$ by the natural group homomorphism mapping $g\in G^{I(X)}$ to the isometry on $Y$ which maps $y$ to $y$ for $y\in Y\setminus I(X)$ and which has the same values as $g$ otherwise.
Let $\nu\in[\mu]$. That means there is an element $g_X\in G^X$ such that $\nu=g_X\mu$. Let $A\subset Y$ be measurable. 
Then for $g_{I(X)}\in G^{I(X)}$ being such that $g_X^{-1}I^{-1}=I^{-1}g_{I(X)}$ it is
\begin{equation}\nonumber
I_\#\nu(A)=I_\# g_X\mu(A)=\mu(g_X^{-1}I^{-1}(A))=\mu(I^{-1}g_{I(X)}(A))=g_{I(X)}^{-1}I_\#\mu(A) \in [I_\#\mu].
\end{equation}
By using Lemma \ref{lem.iso} we now show that $[I_\#]$ indeed gives a distance-preserving map between the Shape spaces: Let $\mu\ne\nu$, then
\begin{eqnarray}\nonumber
D^Y_p([I_\#\mu],[I_\#\nu]) &=& \inf_{g_Y\in G^Y}d^Y(g_Y I_\#\mu,I_\#\nu)=\inf_{g_Y\in G^{I(X)}}d^Y(g_Y I_\# \mu,I_\#\nu) \\ \nonumber
                           &=& \inf_{g_X\in G^X}d^Y(I_\# g_X\mu,I_\#\nu)= \inf_{g_X\in G^X}d^X(g_X\mu,\nu) \\ \nonumber
                           &=& D^X([\mu],[\nu]).
\end{eqnarray}
\end{proof}

\begin{corollary}\label{cor.iso}
$(X,d_X)\cong(Y,d_Y) \Rightarrow (W_p(X)/G^X, D^X_p)\cong (W_p(Y)/G^Y, D^Y_p)$. 
That means, in the given situation the Shape spaces for the full isometry groups are isometric.
\end{corollary}

\begin{proof}
Since every isometry is injective, it only remains to show that $[I_\#]$ is surjective. This follows directly from the surjectivity of the map $I_\#$:
Let $[\eta]\in W_p(X)/G^Y$. Then there is an element $\mu\in W_p(X)$ such that $\eta=I_\#\mu$. That  means, $[\eta]=[I_\#\mu]=[I_\#]([\mu])$.
\end{proof}

\begin{remark}
\begin{itemize}
\item Proposition \ref{prop.iso} still holds true in case either $G^Y$ is substituted by a subgroup which itself contains $G^{I(X)}$ as a subgroup. Or both, $G^X$ and $G^Y$, are substituted by subgroups such that $g\mapsto I\circ g\circ I^{-1}$ is a bijection between these subgroups. For Corollary \ref{cor.iso} the situation can naturally only be changed to the latter scenario.
\item These series of results mean that Shape spaces only depend on the metric structure of the underlying space $X$ and the chosen subgroup of the isometry group. 
\end{itemize}
\end{remark}
\subsection{Geodesic properties of $\Sh_p(X)$}\label{subsec.geod}
As we know from the literature, Wasserstein spaces have an interesting geodesic structure. Naturally, the question arises in which respects the geodesic structure of the Shape spaces is related. 

In this section, let $(X,d)$, $G\subset ISO(X)$ be non-trivial and the action of $G$ on $X$ be such that Theorem \ref{min} holds, i.e. such that $D_p$ is a metric distance. When speaking of geodesics, we always mean constant speed geodesics (as in Definition \ref{df.geod}). 

\begin{prop}\label{prop.shapegeod}
	Let $(\mu_t)_{t\in[0,1]}$ be a geodesic in $W_p(X)$ s.th.	 $D_p([\mu_0],[\mu_1])=W_p(\mu_0,\mu_1)\neq 0$.
    Then $([\mu_t])_{t\in [0,1]}$ is a geodesic as well in $\Sh_p(X)$.
\end{prop}
\begin{proof}
	For every $t,s\in [0,1]$ we have 
	$$D_p([\mu_t],[\mu_s])\leq W_p(\mu_t,\mu_s)=|t-s|D_p([\mu_0],[\mu_1])$$
	which is equivalent to the definition of constant speed geodesics. 
\end{proof}
To be able to formulate a more general statement, we give the following lemma.
\begin{lemma}\label{lem.coef}
Let $(\mu_t)_{t\in[0,1]}$ be any curve in $W_p(X)$. Since in general it holds $D_p([\mu_t],[\mu_s])\leq W_p(\mu_t,\mu_s)$, 
for every $t\neq s\in[0,1]$ there exist $C^{ts}\in\R_{>0}$ such that $C^{ts}D_p([\mu_t],[\mu_s])=W_p(\mu_t,\mu_s)$ whenever $[\mu_t]\neq[\mu_s]$. \hfill $\qed$ 
\end{lemma}

\begin{dfn}[\emph{\textbf{Quotient coefficients}}]
	We call the numbers $C^{ts}$ from Lemma \ref{lem.coef} \emph{quotient coefficients} of the curve $(\mu_t)_{t\in[0,1]}$.
\end{dfn}

\begin{prop}\label{prop.quotiengeod}
	Let $(\mu_t)_{t\in[0,1]}\subset W_p(X)$ be a geodesic. Then $([\mu_t])_{t\in[0,1]}$ is a geodesic, too, if and only if  $[\mu_t]\neq[\mu_s]\ \forall t,s\in[0,1]$ and the quotient coefficients $C^{ts}$ satisfy $C^{ts}=C^{01}\ \forall t\neq s\in[0,1]$.
\end{prop}
Recall that in Definition \ref{df.geod} we demanded $\gamma_0\neq\gamma_1$ for all metric geodesics $(\gamma_t)_{t\in[0,1]}$.

\begin{proof}
If $([\mu_t])_{t\in[0,1]}$ is a geodesic, then by definition $[\mu_t]\neq[\mu_s]\ \forall t\neq s\in[0,1]$. Further, it is
\begin{eqnarray}\nonumber
C^{ts}D_p([\mu_t],[\mu_s]) &=& W_p(\mu_t,\mu_s)=|t-s|W_p(\mu_0,\mu_1)=|t-s|C^{01}D_p([\mu_0],[\mu_1]) \\\nonumber
                       &=& C^{01}D_p([\mu_t],[\mu_s])
\end{eqnarray}
and therefore, $C^{ts}=C^{01}\ \forall t\neq s\in[0,1]$.

Now let  $[\mu_t]\neq[\mu_s]\ \forall t\neq s\in[0,1]$ and $C^{ts}=C^{01}\ \forall t\neq s\in[0,1]$. Then 
\begin{eqnarray}\nonumber
D_p([\mu_t],[\mu_s]) & = & \frac{1}{C^{ts}}W_p(\mu_t,\mu_s)=\frac{1}{C^{01}}W_p(\mu_t,\mu_s)=\frac{1}{C^{01}}|t-s|W_p(\mu_0,\mu_1)\\\nonumber
                     & = & |t-s| D_p([\mu_0],[\mu_1]).
\end{eqnarray}
\end{proof}

\begin{example}
	Consider the function $f:\R_{>0}\rightarrow \R$, $x\mapsto x$ and the revolution of its graph around the $z$-axis. With respect to the induced Euclidean metric structure, the symmetry group consists of the rotations around the $z$-axis. Let $a=(a^x,a^y,a^z)$ be an arbitrary point on the surface of revolution and $b=(b^x,b^y,b^z)$ another point on the surface with $b^y\neq a^y$ and $b^z\neq a^z$. Let $\delta_{a_t}$ be the metric geodesic in the Wasserstein space of this surface between $\delta_a$ and $\delta_b$. Then $D_p(\delta_{a_t},\delta_{a_s})=\sqrt{(a_t^x-a_s^x)^2+(a_t^z-a_s^z)^2}$, while the respective Wasserstein distances also depend on the $y$-components. That means, in this case the quotient coefficients are not equal and thus the projected curve is not a geodesic in Shape space, according to Proposition \ref{prop.quotiengeod}. However, if one would choose $b$ such that $b^z\neq a^z$ but $b^y=a^y$ and $sign(b^x)=sign(a^x)$, the projection of the geodesics between those two points would be a geodesic again in Shape space.
\end{example}

One might wonder whether in case $\mu_1=g_\#\mu_0$, geodesics between $\mu_0$ and $\mu_1$ have the form $\mu_t=g_t\mu_0$, for $g,g_t\in G$. Then, of course, the projection of this curve would only consist of one point, $[\mu_t]=[\mu_0]$. However, this is not true in general. 

\begin{example}
	Consider the graph of the function $f(x)=x^2$ and let it revolve around the $z$-axis. On this manifold with the induced Euclidean metric, the isometry group consists of reflection across and rotation around the $z$-axis. Now consider the Dirac measures $\delta_{a}$, $a=(a_x,a_y,a_z)$ defined by $a_x=-1$ and $a_y=1$ and $\delta_b$, $b=(b_x,b_y,b_z)$ defined by $b_x=1$ and $b_y=1$. $a$ and $b$ are the reflections of one another. A continuous curve connecting $\delta_a$ and $\delta_b$ which is generated by isometries, i.e. a curve of the form  $\delta_t=g_t\delta_a$, can only go along the semicircle in the hyperplane $y=1$, connecting $a$ and $b$. This arc has length $\pi$. But going along the graph of $f(x)$ would be shorter: $\int_{-1}^{1}\sqrt{1+4x^2}dx<\pi$. So $\delta_t=g_t\delta_a$ cannot be a geodesic.
\end{example}

We so far have only asked, when projections of geodesics are geodesics again. Are there, however, other continuous curves which can be projected onto geodesics in Shape space?

\begin{lemma}\label{lem.contihash}
	Let $\mu_t:[0,1]\rightarrow W_p(X)$ be a continuous curve in $W_p(X)$ and $g_t:[0,1]\rightarrow G$ be a continuous curve in $G$. Then $g_{t\#}\mu_t:[0,1]\rightarrow W_p(X)$ is a continuous curve in $W_p(X)$.
\end{lemma}

\begin{proof}
	Let $(a_n)_{n\in\N}$ be a sequence in $[0,1]$ that converges to $a\in [0,1]$. Then $(g_{a_n},\mu_{a_n})$ is a net in $G\times W_p(X)$ which converges to $(g_a,\mu_a)$, since $g_t$ and $\mu_t$ are continuous. 
	We know from Lemma \ref{cont} that the map $G\times W_p(X)\rightarrow W_p(X),\ (g,\mu)\mapsto g_\#\mu$ is continuous in our setting. 
	From this we can deduce that $g_{a_n\#}\mu_{a_n}$ converges to $g_{a\#}\mu_a$, which proves our claim.
\end{proof}

\begin{corollary}\label{cor.projgeod}
Let $X$ be such that $W_p(X)$ is a non-branching geodesic space and such that the path-connected component of the identity element $Id_G$ of $G$ is not trivial. Then there exists a continuous curve $\tilde{\mu}_t:[0,1]\rightarrow W_p(X)$ which is not a geodesic but its projection onto Shape space is.
\end{corollary}

Recall from Subsection \ref{subsec.Wass} that $W_p(X)$ is a non-branching geodesic space as soon as $(X,d)$ is, for $p\in(1,\infty)$.

\begin{proof}
	Since the path-connected component $A$ of $Id_G$ is not trivial, there is an element $\mu_1\in W_p(X)$ such that for $g\in A$, $g\mu_1\neq\mu_1$ (for example $\mu_1=\delta_x$ with $g(x)\neq x$). Now choose $\mu_0\in W_p(X)$ such that $D_p([\mu_0],[\mu_1])=W_p(\mu_0,\mu_1)$ (this is possible due to Theorem \ref{min}). Since $W_p(X)$ is geodesic, there exists a geodesic $\mu_t$ between $\mu_0$ and $\mu_1$ and because of Proposition \ref{prop.shapegeod}, $[\mu_t]$ is a geodesic, too. 
	
	Now let $\gamma_t:[0,1]\rightarrow G$ be a continuous curve from $Id_G$ to $g$, i.e. such that $\gamma_0=Id_G$ and $\gamma_1=g$.
	Let $g_t:[0,1]\rightarrow G$ be such that $g_t=Id_G$ for $0\leq t\leq 1/2$ and $g_t=\gamma_{2t-1}$ for $1/2\leq t\leq 1$.
	From Lemma \ref{lem.contihash} we know that $g_{t\#}\mu_t$ is continuous, but it cannot be a geodesic since $W_p(X)$ is a non-branching space and we chose that $g\mu_1\neq\mu_1$.
\end{proof}

Above, the procedure we used was the pointwise transformation of a geodesic by isometries such that it branches continuously eventually. Another way to construct a curve that is projected onto a geodesic but which itself is none is detailed in the following. It will in particular allow to drop the assumption on the isometry group.

\begin{lemma}\label{lem.projgeod}
 	Let $X$ be such that $W_p(X)$ is a non-branching geodesic space. In $W_p(X)$ let $\mu_t$ be a geodesic such that there is a $g\in G$ with $g\mu_1\neq\mu_1$. Then, if $F:[0,1]\rightarrow[0,1]$ is a continuous function with $F(t)=1$ for all $t\leq t_0<1$ but $F(t)\neq1$ for all $t>t_0$, the mixing $\tilde{\mu}_t:=\left(F(t)+(1-F(t)g)\right)_\#\mu_t$ is not a geodesic, but its projection onto Shape space is.
\end{lemma}

\begin{proof}
	As long as $F(t)=1$, $\tilde{\mu}_t=\mu_t$, but $\tilde{\mu}_t$ branches at $t=t_0$ and can therefore not be a geodesic anymore. But it is $[\tilde{\mu}_t]=F(t)[\mu_t]+(1-F(t))[\mu_t]=[\mu_t]$ and therefore a geodesic by assumption.
\end{proof}

Next to questions regarding the correspondence of curves in Wasserstein space to geodesics in Shape space, it is also interesting to notice whether in Shape space there is always a geodesic between two given points.

\begin{corollary}\label{cor.geodesicshapespace}
	If $W_p(X)$ is geodesic, then $\Sh_p(X)$ is geodesic, too.
\end{corollary}

\begin{proof}
	We argue similarly as in the proof of \ref{cor.projgeod}. 
	Let $\bar{\mu},\bar{\nu}\in\Sh_p(X)$ and $\mu\in W_p(X)$ such that $[\mu]=\bar{\mu}$. Then according to Theorem \ref{min}, there exists a $\nu\in W_p(X)$ such that $[\nu]=\bar{\nu}$ and $W_p(\mu,\nu)=D_p([\mu],[\nu])=D_p(\bar{\mu},\bar{\nu})$. Since we assumed $W_p(X)$ to be geodesic, there exists a geodesic $\mu_t$ in $W_p(X)$ with $\mu_0=\mu$ and $\mu_1=\nu$. Now, according to Proposition \ref{prop.shapegeod}, $[\mu_t]$ is a geodesic in $\Sh_p(X)$ connecting $\bar{\mu}$ and $\bar{\nu}$.
\end{proof}

We know from the literature that $W_2(X)$ is geodesic whenever $(X,d)$ is geodesic.

The geodesic we have found in the proof of \ref{cor.geodesicshapespace} between $\bar{\mu}$ and $\bar{\nu}$ need not be unique. A priori it is not clear whether geodesics between other elements of the equivalence classes $\bar{\mu}$ and $\bar{\nu}$ which have Wasserstein distance $D_p(\bar{\mu},\bar{\nu})$, are projected onto the same curve in Shape space. And, quite in general, it can still be that there are geodesics in Shape space which cannot be written as the projection of a geodesic, or any other continuous curve, in Wasserstein space at all.
\section{Tangent spaces on $\Sh_2(M)$} \label{ch.infSh}
%

The goal of this section is to construct a tangent space on Shape space $\Sh(M)$, where $(M,h)$ is a smooth, connected, complete Riemannian manifold. 
Before doing this, it is instructive to have a look at the case of finite dimensional manifolds, where the situation is already well explored:
In case a finite dimensional Lie group $G$ acts smoothly, freely and properly on a finite dimensional smooth manifold $M$, the quotient space $M/G$ is a topological manifold of dimension dim $M/G=$ dim $M\ -$ dim $G$. The quotient space furthermore has a unique smooth structure with respect to which the quotient map $\pi:M\rightarrow M/G$ is a smooth submersion (see \cite{Lee2001}). 
The orbits $G.x$, $x\in M$, then, can proven to be embedded submanifolds of $M$ and the kernel of the linear and surjective map $d\pi_x:T_xM\rightarrow T_{\pi(x)}\left(M/G\right)$ is precisely $T_x(G.x)$ (for these statements see again \cite{Lee2001}, Lemma 5.29 and Theorem 7.5). This yields the following isomorphism of linear spaces: 
\begin{equation}\label{eq.finitequotienttangent}
T_{[x]}\left(\bigslant{M}{G}\right)\cong\ \bigslant{T_xM}{T_x\left(G.x\right)}.
\end{equation}
This means, qualitatively, in $T_{[x]}(M/G)$ all directions showing ``orbitwards'' are modded out.

Of course, we cannot apply this statement to our case - $W(M)$ is an infinite dimensional, only \emph{formal}, Riemannian manifold and also our $G$-action is not free in general\footnote{Just consider  a normal distribution on $\R^3$ centered around the origin, which is invariant under rotations.}. However, we can  mod out \emph{by hand} those directions of $T_\mu W(M)$ that point towards the orbit of $\mu$. 
\newline

Recall that the group of isometries $ISO(M)$ of a finite dimensional Riemannian manifold $M$ is a Lie group, due to a theorem by Myers and Steenrod, \cite{Myers_1939}. And in case $M=\R^n$, $ISO(M)=O(n)\ltimes\R^n$, where $O(n)$ is the orthogonal group of $\R^n$ and $\ltimes$ denotes the semidirect product of groups. 
We will find the above mentioned orbitward directions on Wasserstein space to be induced by the Lie algebra $\mathfrak{iso}(n)$ of $ISO(\R^n)$. One can show that $\mathfrak{iso}(n)=\mathfrak{so}(n)\ltimes \R^n$, where $\ltimes$ denotes the semidirect sum of algebras. (See e.g. \cite{Fecko2006} for a reference on this result and for a definition of the semidirect sum. However, we will be more explicit, below.) By $\mathfrak{so}(n):=\{M\in \R^{n\times n}\mid M^t=-M\}$ we denote the vector space of skew-symmetric matrices with the matrix-commutator $[M,N]:=M\cdot N-N\cdot M$.
It is $\dim \mathfrak{iso}(n) = \dim ISO(\R^n)=\frac{n(n+1)}{2}$.
\newline 

We assume that the reader is familiar with the basic notions of Lie theory.\footnote{For more background knowledge on Lie theory consult for example \cite{Baum2009} (in German) or \cite{Lee2001}.} However, we would like to recall the definition of fundamental vector fields, as they will be important for us.

Let $\mathfrak{g}$ be the Lie algebra of a Lie group $G$, $X\in\mathfrak{g}$ a left-invariant vector field 
and $\gamma_X(t)$ the maximal integral curve of $X$ starting at the identity element, i.e. with $\gamma_X(0)=e$. Then $\varphi_X(t)$ is defined for all values $t\in\R$ and the map $\varphi_X:\R\rightarrow G$ is a homomorphism of Lie groups, i.e. $\gamma_X(0)=e$, $\varphi_X(s+t)=\varphi_X(s)\cdot \varphi_X(t)$ and $\varphi_{sX}(t)=\varphi_X(s\cdot t)$ for all $s,t\in\R$.

\begin{dfn}[\emph{\textbf{Exponential map}}]
  The map 
  \begin{eqnarray}\nonumber
  \exp:\mathfrak{g} &\longrightarrow& G\\\nonumber
            X  &\longmapsto& \varphi_X(1)
  \end{eqnarray}
  is called the \emph{exponential map} of the Lie group $G$. 
  
   The curve $\R\ni t\mapsto\exp(tX)=\varphi_{tX}(1)=\varphi_X(t)$, then, is the maximal integral curve of $X$ through $e$ and the curve $\R\ni t\mapsto g\cdot\exp(tX)$ the maximal integral curve of $X$ through $g\in G$.
\end{dfn}

We now assume that the Lie group $G$ acts on a smooth manifold $M$. Then, every element of the Lie algebra $\mathfrak{g}$ of $G$ defines a certain vector field on $M$, the so called \emph{fundamental vector field}.

\begin{dfn}[\emph{\textbf{Fundamental vector field}}]
	Let $G$ act on $M$ on the left and let further be $X\in\mathfrak{g}$. The vector field $\tilde{X}\in\Gamma^{\infty}(M)$ on $M$ defined by
	\begin{equation}\nonumber
	\tilde{X}(x):=\frac{d}{dt}\left(\exp(-tX)\cdot x\right)|_{t=0}
	\end{equation}
	is called the \emph{fundamental vector field} corresponding to $X$.\footnote{There is a similar definition for right-actions of $G$. But here, we will only be interested in left-actions.}
\end{dfn}
The map $\mathfrak{g}\rightarrow\Gamma^{\infty}(M),\ X\mapsto \tilde{X}$ is linear and $\widetilde{[X,Y]}=[\tilde{X},\tilde{Y}]$, where the latter is the commutator for vector fields, i.e. the Lie bracket. In particular, the set of all fundamental vector fields form a Lie-subalgebra of the Lie-algebra of all smooth vector fields on $M$ and is Lie-algebra homomorphic to $\mathfrak{g}$.

Is $G$ the isometry group of a complete Riemannian manifold $(M,h)$, the Lie algebra of fundamental vector fields corresponding to $\mathfrak{g}$ is precisely the Lie algebra of Killing vector fields\footnote{A vector field $X$ is a Killing vector field of the Riemannian metric $h$, in case $\mathcal{L}_Xh=0$, where $\mathcal{L}$ denotes the Lie derivative.}, the so called \emph{Killing algebra} on $(M,h)$.
In case $M=\R^n$ and $h$ is the Euclidean metric tensor, one can calculate the general Killing vector field (\cite{Fecko2006}), and thus the general fundamental vector field, to be of the form
\begin{equation}\label{eq.formkilling}
    \xi^{(A,a)}=\frac{1}{2}(Ah)^{ij}M_{ij}+a^iP_i,
\end{equation}
where $A\in\mathfrak{so}(n)$ and $a\in\R^n$. Here, we have used the Einstein summation convention. The vector fields 
\begin{equation}
M_{ij}=-M_{ji}= x_i\partial_j-x_j\partial_i\ \text{ and }\ P_i=\partial_i
\end{equation}
constitute a basis of the Killing algebra.
\newline 

Recall from Section \ref{sec.dyn}, in particular Theorem \ref{them.ac}, that the tangent space on $W(M)$ is constructed using the weak continuity equation. A family of vector fields $v_t$ is considered to be tangent along an absolutely continuous curve $\mu_t$ if the weak continuity equation holds and if $\|v_t\|_{L^2(\mu_t)}$ is minimal among all the possible solutions. 
Now, to identify all those tangent vectors that show in the direction of the orbit, let us consider the pushforward of a measure $\mu$ by the flow of the left-invariant vector fields on $G=ISO(M)$.

\begin{prop}\label{prop.fundconti}
	Let $X\in\mathfrak{g}$ and $\tilde{X}$ be the corresponding fundamental vector field on $M$. Further, let $\mu\in W(M)$ and $\mu_t:=\exp(-tX)_{\#}\mu$, $0\leq t\leq 1$. Then, the tuple $(\mu_t,\tilde{X})$ is a solution of the weak continuity equation \eqref{eq.wcont}, i.e. 
	\begin{equation}\nonumber
	\int_0^1\int_M\left(\frac{\partial}{\partial t}\varphi(x,t)+h(\nabla\varphi(x,t),\tilde{X}(x))\right)\ d(\exp(-tX)_{\#}\mu)(x)dt =0.
	\end{equation} 
\end{prop}

\begin{proof}
	The integral curve $\exp(-tX)$ on $G$ defines the flow $\exp(-tX):M\rightarrow M$ on $M$ which is a diffeomorphism for every $t$ and $\exp(0)=Id$. It is in particular the flow of the fundamental vector field $\tilde{X}$, so that the statement of the proposition follows from Theorem \ref{thm.conti}.
\end{proof}

Let us note that because the tangent map of the orbit map $g_x:G\rightarrow G.x,\ g\mapsto g.x$ at point $e$, i.e. $(dg_x)_e:T_eG\rightarrow T_xG.x$, is surjective, every curve that is initially tangential to $G.x$ infinitesimally coincides with the integral line of a fundamental vector field starting at $x$. This is why in our case we will settle with those kinds of curves considered in Proposition \ref{prop.fundconti}.
By doing this, we naturally only consider the path-connected component of $Id_G$ with respect to which directions it is going to be modded out by.
\newline

The next question to be asked is whether $\mu_t$ from Proposition \ref{prop.fundconti} is indeed absolutely continuous in $W(M)$. According to Theorem \ref{them.ac} this is so, in case $\int_0^1 \|\tilde{X}\|_{L^2(\mu_t)}\ dt < \infty$. 

We will give a positive answer for this only in case $(M,h)$ is the Euclidean space. However, we conjecture that the same holds true also for compact manifolds and Riemannian manifolds with nonnegative curvature. In the hyperbolic case a restriction to the subspace of compactly supported measures seems necessary.

\begin{prop}\label{prop.regfund}
	Let $\tilde{X}$ be a fundamental vector field corresponding to an element $X\in\mathfrak{iso}(n)$ and let $\mu\in W(\R^n)$. Then 
	$$\|\tilde{X}\|_{L^2(\mu)}<\infty.$$
\end{prop}

\begin{proof}
	First we point out that since $\mu\in W(M)$, it is $\int_{\R^n}\|x\|^2\ d\mu(x)<\infty$. Because of this, also
	$\int_{\R^n}x_i\ d\mu(x)<\infty$ and $\int_{\R^n}x_ix_j\ d\mu(x)<\infty$, where $x_i$ denotes the $i$-th component of $x$. 
	
	In \eqref{eq.formkilling} we have seen that every fundamental vector field $\tilde{X}$ on $\R^n$ is of the form $\tilde{X}=\frac{1}{2}(Ah)^{ij}M_{ij}+a^iP_i$. So, calculating the Euclidean norm of $\tilde{X}$, we obtain for $a=(a_1,...,a_n)$,
	\begin{equation}\nonumber
	\|\tilde{X}\|^2= \|\frac{1}{2}A^{ij}M_{ij}\|^2+\|a\|^2+\langle A^{ij}M_{ij},a^iP_i\rangle.
	\end{equation}
	The integral over the second term is finite. The third term is a linear combination of the components $x_i$, whose integral is again finite. 
	The first term consists of a linear combination of terms of the form $x_i^2+x_j^2$ plus a linear combination of terms of the form $x_ix_j$. The integral over both of these terms is also finite, which means that in total $\int_{\R^n}\|\tilde{X}\|^2d\mu = \|\tilde{X}\|_{L^2(\mu)}^2<\infty$.
\end{proof}

\begin{lemma}\label{lem.norminvariance}
	Let $G$ be the Lie group of isometries of a Riemannian manifold $M$ and $\mathfrak{g}$ its Lie algebra. Further, let $\tilde{X}$ be the fundamental vector field corresponding to $X\in\mathfrak{g}$. If $\mu\in \PM(M)$ and $\mu_t:=(\exp(-tX))_\#\mu$, then, for all $t$,
	$$ \|\tilde{X}\|_{L^2(\mu_t)}^2 = \|\tilde{X}\|^2_{L^2(\mu)}. $$
\end{lemma}

\begin{proof}
		\begin{eqnarray}\nonumber
        \|\tilde{X}\|_{L^2(\mu_t)}^2 &=& \int_{\R^n} \|\tilde{X}_x\|^2\ d((\exp(-tX)_\#\mu)(x) = \int_{\R^n} \|\tilde{X}_{\exp(-tX)y}\|^2\ d\mu(y)\\\nonumber
        &=& \int_{\R^n} \|d(\exp(-tX))\tilde{X}_{y}\|^2\ d\mu(y) = \int_{\R^n} \|\tilde{X}_y\|^2\ d\mu(y)\\\nonumber
        &=& \|\tilde{X}\|^2_{L^2(\mu)}.
\end{eqnarray}	
\end{proof}

\begin{corollary}\label{cor.regfund}
	Let $X\in\mathfrak{iso}(n)$, $\tilde{X}$ the corresponding fundamental vector field, $\mu\in W(M)$ and $\mu_t:=(\exp(-tX))_\#\mu$.
    Then $$\int_0^1 \|\tilde{X}\|_{L^2(\mu_t)}\ dt < \infty.$$
\end{corollary}

\begin{proof}
	Combining Lemma \ref{lem.norminvariance} and Proposition \ref{prop.regfund}, we get 
		\begin{equation}\nonumber
		\int_0^1 \|\tilde{X}\|_{L^2(\mu_t)}\ dt = \int_0^1 \|\tilde{X}\|_{L^2(\mu)}\ dt = \|\tilde{X}\|_{L^2(\mu)} <\infty.
		\end{equation}
\end{proof}

\hfill\break
But although $\tilde{X}$ is an element of $L^2(\mu_t)$ for all $t$ and fulfills the regularity condition of Corollary \ref{cor.regfund}, it is not necessarily the minimal vector field along $\mu_t$. (Consider for example for $n=2,$ $\mu_t:=\exp(-tX)_\#\mu$ where $\mu$ has a normal distribution and $\exp(-tX)$ consists of rotations around the origin. Then $\mu_t=\mu$ for all $t$ and the minimal vector field is $0$.) But we can apply the orthogonal linear projection $P^\mu:L^2(\mu) \rightarrow  T_\mu W(\R^n)$, so that $P^{\mu_t}(\tilde{X})\in T_{\mu_t}W(\R^n)$ and the continuity equation still holds.
This is why we define
$$U_\mu:=\{P^{\mu}(\tilde{X})\mid \tilde{X} \text{ is a fundamental vector field} \}\subset L^2(\mu).$$ 
$U_\mu$ is a complete subspace of $T_\mu W(\R^n)$ and contains all the vectors pointing towards the orbit of $\mu$. With this definition, we can finally propose a notion for $T_{[\mu]}\Sh(\R^n)$.

\begin{dfn}[\emph{\textbf{Tangent space on $\Sh(\R^n)$}}]\label{dfn.tsh}
	We define the tangent space on $\Sh(\R^n)$ at point $[\mu]$ to be the set
	\begin{equation}\nonumber
	T_{[\mu]}\Sh(\R^n):=\ \bigslant{T_\mu W(\R^n)}{U_\mu}.
	\end{equation}
\end{dfn}

\begin{prop}\label{prop.independence}
	It is 
	\begin{equation}\nonumber
	\bigslant{T_\mu W(\R^n)}{U_\mu}\ \cong\ \bigslant{T_{g\mu} W(\R^n)}{U_{g\mu}}.
	\end{equation}
	In this sense, Definition \ref{dfn.tsh} is independent of the choice of the representative.
\end{prop}

\begin{proof}
	We know that $dg:T_\mu W(\R^n)\rightarrow T_{g\mu}W(\R^n)$ is an isomorphism for every $g\in G$. It remains to show that $dg$ also maps $U_\mu$ isomorphically onto $U_{g\mu}$. 
	To show this, let us recall the known fact that if $\tilde{X}$ is a fundamental vector field, then so is $dg(\tilde{X})$ (see for example \cite{Baum2009}, Theorem 1.25). 
	Since for $v=v^\top\oplus v^\perp\in T_\mu W(\R^n)\oplus T^\perp_\mu W(\R^n)$, $dg(v)=dg(v^\top)\oplus dg(v^\perp)$ with $dg(v^\top)\in T_{g\mu} W(\R^n)$ and $dg(v^\perp)\in T^\perp_\mu W(\R^n)$, it is $dg\circ P^\mu= P^{g\mu}\circ dg$.
	So, if $v\in U_\mu$, then there is a fundamental vector field $\tilde{X}$ such that $v=P^\mu(\tilde{X})$ and $dg(v)=dg(P^\mu(\tilde{X}))= P^{g\mu}(dg(\tilde{X}))$ and the statement holds true.
\end{proof}


\begin{example}	
	\begin{equation}\nonumber
	T_{[\delta]}\Sh(\R^3)=\bigslant{\R^3}{U_\mu} \cong \{0\}.
	\end{equation}
	On can interpret this in a sense that there is no non-trivial classical one-particle motion in Shape space.
\end{example}

\begin{remark}
	Definition \ref{dfn.tsh} is meaningful for every $\mu\in W(M)$ on a general Riemannian manifold, as long as $\int_0^1 \|\tilde{X}\|_{L^2(\mu_t)}\ dt < \infty$ for every $\mu_t$ of the form $\mu_t=(\exp(-tX))_\#\mu$. In particular, Proposition \ref{prop.independence} does not use any properties specific to $\R^n$. Just substitute $M$ for $\R^n$.
\end{remark}
\section*{Acknowledgements}
I am indebted to my PhD supervisor Prof. Thomas Schick from the Georg-August-Universität Göttingen for his guidance and mathematical advice throughout my PhD.

\bibliographystyle{alpha}       
\bibliography{da,C:/Users/Admin/Documents/Archives/Bibliothek/Mathe.bib,C:/Users/Admin/Documents/Archives/Bibliothek/Physics.bib}

\end{document}